\documentclass[11pt]{amsart}
\usepackage[colorlinks=true,citecolor=blue,linkcolor=blue, backref=page]{hyperref}
\usepackage{amssymb,verbatim,amscd,amsmath,graphicx}
\usepackage{enumerate}
\usepackage{graphicx}
\usepackage{caption}
\usepackage{subcaption}
\usepackage{pdfsync}
\usepackage{color,colortbl}
\usepackage{fullpage}
\usepackage{bbm}
\usepackage{comment}
\usepackage{multirow}
\usepackage{enumitem}
\usepackage{comment}
\numberwithin{equation}{section}
\newtheorem{introthm*}{Main Result}
\newtheorem{theorem}{Theorem}[section]
\newtheorem{lemma}[theorem]{Lemma}

\newtheorem{corollary}[theorem]{Corollary}

\theoremstyle{definition}
\newtheorem{definition}[theorem]{Definition}
\newtheorem{conjecture}[theorem]{Conjecture}
\newtheorem{def-prop}[theorem]{Definition-Proposition}

\newtheorem{remark}[theorem]{Remark}

\newtheorem*{acknowledgement}{Acknowledgements}

\DeclareMathOperator{\reg}{reg}

\DeclareMathOperator{\Ass}{Ass}

\DeclareMathOperator{\Proj}{Proj}

\renewcommand{\AA}{{\mathbb A}}

\newcommand{\PP}{{\mathbb P}}

\newcommand{\NN}{{\mathbb N}}

\newcommand{\XX}{{\mathbb X}}

\newcommand{\kk}{{\mathbbm k}}

\def\L{{\mathcal L}}

\def\P{{\mathcal P}}

\def\mm{{\frak m}}

\def\pp{{\frak p}}

\def\ahat{\widehat{\alpha}}

\def\1{{\bf 1}}
\def\0{{\bf 0}}

\begin{document}
	
	\title{Demailly's Conjecture for general and very general points}
	
	\author{Sankhaneel Bisui}
	\address{Arizona State University \\ School of Mathematical and Statistical Sciences \\
		Tempe, AZ 85287-1804, USA}
	\email{sankhaneel.bisui@asu.edu} 
	
	\author{Dipendranath Mahato}
	\address{Tulane University, Department of Mathematics,
		6823 St. Charles Avenue
		New Orleans, LA 70118, USA}

	\email{dmahato@tulane.edu}

	\keywords{Demailly's Conjecture, Waldschmidt Constant, Waldschmidt Decompostion Ideals of Points, Symbolic Powers, Containment Problem}
	
	\begin{abstract}
	We prove that at least $\left( \dfrac{(1+\epsilon)2m}{N-1}+1+\epsilon \right)^N$, where $0\leqslant \epsilon <1$, many general points,  satisfy Demailly's conjecture. Previously, it was known to be true for at least $(2m+2)^N$ many general points \cite{bghndemailly}. We also study Demailly's conjecture for $m=3$ for ideal defining general and very general points.  		
\end{abstract}
\maketitle
\section{Introduction}
Let $P_1, \dots, P_s$ be $s$ reduced points in $\PP^N$. Then, the interpolation problem seeks to find the least degree of a hypersurface that passes through the points at least $m$ times. Since homogeneous polynomials define hypersurfaces, it is equivalent to finding the least degree of a homogeneous polynomial that vanishes at those points with multiplicity at least $m$. Let $I$ be the defining ideal of the reduced points 
$P_1, \dots, P_s$ then by Zariski-Nagata theorem, the $m$-th symbolic power of $I$, denoted by $I^{(m)}$ is the collection of all homogeneous polynomials that vanishes at those points with multiplicity $m$. 
It has been proved that the Waldschmidt constant, denoted by $ \displaystyle \ahat(I)=\lim_{m \to \infty} \dfrac{\alpha(I^{(m)})}{m}=\inf _{m \to \infty} \dfrac{\alpha(I^{(m)})}{m} $ exists. Chudnovsky \cite{chudnovsky1979singular} proposed a conjectural bound to the interpolation problem, which Demailly\cite{demailly1982formules} generalized to the following statement about the lower bound of the Waldschmidt constant of the defining ideal of the points. 
\begin{conjecture}\label{conjecture: Demailly}
	If $I$ is the defining ideal of $s$ reduced points in $\PP^N$ then $$\ahat(I)\geqslant \dfrac{\alpha(I^{(m)})+N-1}{m+N-1}, \text{ for all } m \geqslant 1. $$
\end{conjecture}   
For $m=1$ Conjecture \ref{conjecture: Demailly} is Chudnovsky's conjecture and is very well studied, e.g., see,  \cite{dumnicki2015containments,FMXchudnovsky, dtChudnovsky2017containment,bisuinguyen2021chudnovsky, bisuinguyen2021chudnovsky}. Demailly's conjecture is also well-studied. It has been proved for points in $\PP^2$ by Esnault and Viehweg \cite{esnaultviehweg}. Using the ideal containment method, Bocci, Cooper, and Harbourne\cite{bocci2014containment} proved that a square number of points in $\PP^2$ satisfies Demailly's conjecture. Malara, Szemberg, and Szpond \cite{malara2018conjecture} proved that at least $(m+1)^N$ very general points satisfy Demailly's conjecture which Dumnicki, Szemberg, and Szpond \cite{dumnicki2024waldschmidt} improved for at least $m^N$ many very general points. Chang and Jow \cite{chang2020demailly} proved that for a fixed $m$ the ideal defining $s\geqslant (\dfrac{2\epsilon(m-1)}{N-1} +2)^N$ very general points, where $0\leqslant \epsilon <1$  satisfies Demailly's conjecture. Later, Bisui, Grifo, Hà, and Nguy$\tilde{\text{\^e}}$n \cite{bghndemailly} proved that for a fixed $m$ if $I$ is the defining ideal of a general set of at least $(2m+2)^N$  points in $\PP^N_\kk$ then $I$ satisfies Conjecture \ref{conjecture: Demailly}. It was proved that points forming star configurations \cite{bghndemailly,malara2018conjecture} and Fermat configuration \cite{Ng23b} satisfy Demailly's conjecture. Demailly-like bounds were proved for determinantal ideals \cite{bghndemailly}, and some special configurations of lines \cite{nguyen2023initial}. 

A property is said to hold for $s$ very general points if it holds for points in a countable infinite intersection of open Zariski sets in the Chow variety of $s$ points. A property is said to be true for $s$ general points if it holds for points only in one Zariski open set. We refer to \cite[Page 4]{FMXchudnovsky} for more details.  

Chang and Jow \cite{chang2020demailly} proved that given an $m \geqslant 1$, there is a countable infinite intersection of open sets in the Chow variety of $s$ points, where $s\geqslant (\dfrac{2\epsilon(m-1)}{N-1} +2)^N$ and $0\leqslant \epsilon <1 $, such that the points from that intersection satisfies Demailly's conjecture. In this manuscript, we prove that given $m\geqslant 1$, there is an open set in the Chow variety of $s$ points such that the points from that open set satisfy Demailly's conjecture.  
The following theorem gives the bound on the number of points. Since we go from a countable infinite intersection of open sets to only one open set, we have more restrictions. 
\begin{introthm*}[Theorem \ref{theorm: Demailly general}]
	Let $\kk$ be an algebraically closed field and $N\geqslant 3$. Fix an integer $m\geqslant 1$. If $I$ denotes the ideal defining a set of at least $\left( \dfrac{(1+\epsilon)2m}{N-1}+1+\epsilon \right)^N $ general points in $\PP^N_\kk$, where $0 \leqslant \epsilon <1$, then
	$$\ahat(I)\geqslant\dfrac{\alpha(I^{(m)})+N-1}{m+N-1}.$$
\end{introthm*} 

As mentioned in Remark \ref{remark: Demailly general what has been improved}, Theorem \ref{theorm: Demailly general} implies that for a given $m$, the defining ideal $I$ of  at least $4^N$ general points in $\PP^N, N\geqslant 2m+1$ satisfies $\ahat(I)\geqslant \dfrac{\alpha(I^{(m)})+N-1}{m+N-1}$. This improves the restriction in \cite[Theorem 2.9]{bghndemailly} on the number of points in $\PP^N$, where $N$ is large. 
\\

In \cite{bn2024lower}, a specific case of Conjecture \ref{conjecture: Demailly}, namely, Demailly's conjecture for $m=2$, has been studied. Inspired by that, in this manuscript, we concentrate on the following particular case of Conjecture \ref{conjecture: Demailly}, namely, Demailly's conjecture for $m=3$.  
\begin{conjecture}\label{conjecture: Demailly3}
	If $I$ is the defining ideal of $s$ reduced points in $\PP^N$ then $$\ahat(I)\geqslant \dfrac{\alpha(I^{(3)})+N-1}{N+2}. $$
\end{conjecture}

Conjecture \ref{conjecture: Demailly3} is known for the following cases:  
at least $3^N$ very general points \cite{dumnicki2024waldschmidt},  and at most $N+3$ linear general points \cite{nagel2019interpolation} in $\PP^N$. In $\PP^3$, it is known for at least $8^3$ general points \cite{bghndemailly}, in $\PP^4$, it is known for at least $7^4$ many general points \cite{bghndemailly} and in $\PP^N$ it is known for at least $6^N$ general points when $N\geqslant 5$ \cite{bghndemailly}. Due to the difficulty of obtaining the initial degree of the third symbolic power and the Waldschmidt constant, we must investigate each case thoroughly. We improve the previous results as follows.

\begin{introthm*} We prove Conjecture \ref{conjecture: Demailly3} holds for the following cases:
	\begin{enumerate}[label=(\roman*)]
		\item  (Corollary \ref{corollary: Demailly m=3 from demailly general}) $s$ general points in $\PP^N$, where
		\begin{enumerate}[label=(\alph*)]
			\item $s\geqslant 6^4$, where $N=4$
			\item $s \geqslant 5^N$, where $N=5$, or $6$
			\item $s\geqslant 4^N$, where $N \geqslant 7$.
		\end{enumerate}
		\item (Corollary \ref{corollary: demailly m=3 very general in lower PN}) Any number of very general points in $\PP^N$, where $3 \leqslant N\leqslant 7$. 
		\item (Remark \ref{remark: demailly for few cases}) at least $\left\lceil \dfrac{{N+\ell \choose N}}{{ N+2 \choose N}}   \right\rceil $ very general points in $\PP^N$, where $N \geqslant 8$, and $\ell= \dfrac{5+\sqrt{8N+17}}{2}$. 
	\end{enumerate}
\end{introthm*}
We also proved that if we consider a slightly weak bound, every very general point in projective spaces will satisfy it. 
\begin{introthm*}[Theorem \ref{thm: demailly very general when m=3}]
	If $I$ is the defining ideal of $s$ very general points in $\PP^N$, then $$\ahat(I)\geqslant \dfrac{\alpha(I^{(3)})+N-2}{N+2}. $$	
\end{introthm*}
Now, we briefly describe our strategy. 
The main tool in \cite{bghndemailly} was ideal containment. Demailly's conjecture follows from the following containment conjecture \cite{HarbourneHuneke}: $$I^{(r(m+N-1))} \subset \mm ^{r(h-1)}\left(I^{(m)}\right)^r, \text{ for all } r,m\geqslant 1. $$
It has been shown in \cite[Lemma 2.2]{bghndemailly} that a stronger containment namely \[I^{(c(m+N-1)-N+1)} \subset \mm ^{r(N-1)}\left(I^{(m)}\right)^r \tag{S}\label{S},\] for a fixed positive integer $c$ and for all $m\geqslant 1$, will imply Demailly's Conjecture \ref{conjecture: Demailly}. 
We proved that if $s\geqslant \left( \dfrac{(1+\epsilon)2m}{N-1}+1+\epsilon \right)^N$, where $0\leqslant \epsilon<1$, then the ideal $I$ defining $s$ generic points satisfies the above stronger containment (\ref{S}) for some $c$. To get the containment, we showed that the appropriate initial degree bounds the regularity of $m$-th symbolic power. To bound the regularity of the $m$-th symbolic powers, we used the well-known Lemma \ref{lemma: regularity bound trung}, which says that 
$$ \text{ if } (s-1){N+m-1\choose N} < {w+N\choose N} \text{ then } \reg(I^{(m)})\leqslant m+w  .$$

To prove Demailly's conjecture concerning the initial degree of the third symbolic power, our strategy is to bound the initial degree of the third symbolic powers as follows: 
$$ \text{ if }{d+N-1\choose N}\leqslant  s{N+2\choose N} < {d+N\choose N} \text{ then } \alpha(I^{(3)})\leqslant d.  $$ 

Then, we break the number of points in a balanced way and use Lemma \ref{lemma: Waldschmidtdecomp} to obtain the proper value of the Waldschmidt constant. This method, together with induction, gives the desired results. The technique has the potential to be applied for any $m$. But one has to be cautious for large $m$ and small number of points.

This article is outlined as follows. In Section \ref{secction: preliminaries}, we recall the significant results that we used. In Section \ref{section: Demailly's conjecture for general points}, we prove Demailly's conjecture for a large number of general points, improving the number of points. In Section \ref{section: compairing waldschmidt constant with the third symbolic powers}, we compared the Waldschmidt constant with the initial degree of the third symbolic powers of the defining ideal of very general points. In the last section, Section \ref{section: computations in lower dimensions}, we showed specific computations for points in lower dimensional projective spaces $\PP^N$, where $3 \leqslant N \leqslant 7$.    

\begin{acknowledgement}
	The authors are thankful to Huy Tài Hà and Thái Thành Nguy$\tilde{\text{\^e}}$n for various conversations. The authors are grateful to Huy  Tài Hà for valuable comments on the draft of the paper. This is part of the second author's thesis. 
\end{acknowledgement}

\section{Preliminaries}\label{secction: preliminaries}	  
In this section, we recall some important results and terminologies that we will use. We start by defining the symbolic powers of ideals in a polynomial ring. 
\begin{definition}
	If $I \subset R=\kk[\PP^N] $ is homogeneous ideal then the $m$-th symbolic power of $I$ is defined by $$I^{(m)}=\cap_{\pp\in \Ass(I)} I^m R_\pp\cap R.$$
	The Waldschmidt constant of the ideal $I$ is defined by $\displaystyle \ahat(I) =\lim_{m \to \infty}\dfrac{\alpha(I^{(m)})}{m}=\inf \dfrac{\alpha(I^{(m)})}{m}$, where $\alpha(I^{(m)})$ is the initial degree of the ideal $I^{(m)}$.
\end{definition} 

We briefly recall the description of general and very general points. We refer to \cite{FMXchudnovsky} for details. Let $z_{i,j}$ be  $s(N+1)$ many indeterminates and $\kk(z_{i,j})$ be the quotient field of $\kk[z_{i,j}]$.  We consider $s$ points in $\PP^N_{\kk(z_{i,j})}$, namely, $P_1(\underline{z})=[z_{1,0}, \dots,z_{1,N} ], \dots, P_s(\underline{z})=[z_{s,0}, \dots,z_{s,N}  ]$.  Let  $\XX(\underline{z})=\{P_1(\underline{z}), \dots,  P_s(\underline{z}) \}. $ Let $\underline{a}\in \AA^{s(N+1)}_\kk$ be any point. Then by coordinate-wise substitution we get $\XX(\underline{a})=\{P_1(\underline{a}), \dots, P_s(\underline{a}) \}$. The set $\XX(\underline{z})$ is said to be the set of $s$ generic points in $\PP^N_{\kk(\underline{z})}$ and $\XX(\underline{a})$ is said to be the set of points obtained via specialization. A property $\P$ holds for $s$ very general (or general) points if there is a countable infinite intersection of open sets (or one open ) in $\AA^{s(N+1)}$  such that points obtained from the set after specialization satisfies the property. 
We refer \cite{bghnchudnovsky} for detailed work on the behavior of the symbolic power of ideals defining points with specialization. 

Now, we recall one technical lemma that we will use. 
\begin{lemma}[Theorem 4.1 \cite{dumnicki2024waldschmidt}]\label{lemma: Waldschmidtdecomp}
	Denote $\ahat(\PP^N, r)$ the Waldschmidt constant of the ideal of $r$ very general points in $\PP^N$. Let $N\geqslant 2$ and $k\geqslant 1$. Assume that for some integers $r_1,\ldots ,r_{k+1}$ and rational numbers $a_1,\ldots ,a_{k+1}$ we have 
	$$\ahat(\PP^{N-1}, r_j) \geqslant a_j, \text{ for all } j=1,\ldots ,k+1,$$
	$$k\leqslant a_j \leqslant k+1, \text{ for } j=1,\ldots ,k, \ \  a_1>k, \ \  a_{k+1}\leqslant k+1.$$
	Then, $$\ahat(\PP^{N}, r_1+\ldots +r_{k+1}) \geqslant \left(1 - \sum_{j=1}^k \dfrac{1}{a_j} \right)a_{k+1} +k.$$
\end{lemma} 

Let $P_1, \dots, P_s$ be $s$ many reduced generic points in $\PP^N$ and $\pp$ be the ideal defining the point $P$. The scheme defined by the ideal $\pp_1^{m_1}\cap \dots \cap\pp_s^{m_s}$ is said to be the fat point scheme, and it is denoted by $Y=m_1P_1+\dots +m_sP_s$. Let $\L(d)=\Proj\left([\pp_1^{m_1}\cap \dots \cap\pp_s^{m_s}]_d\right)$.  
It is said that the $Y$ has \emph{good postulation} in degree $d$ if the following conditions are satisfied: 
\begin{itemize}
	\item if $\deg(Y) \leqslant {N+d \choose N}$ then $h^1\left(\PP^N, \L(d) \right)=0$. 
	\item  if $\deg (Y) \geqslant {N+d \choose d}$ then $h^0\left(\PP^N, \L(d) \right)=0$. 
\end{itemize}
If the fat point scheme $mP_1+\dots+mP_s$ is in good postulation in degree $d$, then the regularity of the symbolic powers of the defining ideal of $P_1, \dots, P_s$ can be easily obtained via binomial inequalities.

$$  {d+N-1\choose N} \leqslant s{N+m-1\choose N} < {d+N\choose N} \text{ then } \reg(I^{(m)})\leqslant d+1  .$$

Generally, proving that fat point schemes are in good postulation in degree $d$ is difficult. 
Alexander-Hirscowitch \cite{alexander1988singularites,alexander1992lemme,alexander1992methode,alexander1995polynomial} theorem showed that $2P_1+\dots +2P_s$ are in good postulation except few cases. It has been proved that quintuple points in $\PP^3$ are in good postulation in degree $d$ except in a few cases \cite{ballico2009postulation, ballico2012postulation}. 

The following is the other useful bound for the regularity of the $m$-th symbolic power defining the ideal of points. 
\begin{lemma}\label{lemma: regularity bound trung}\cite[Theorem 2.4]{TrungValla}
	Assume that $\kk$ is an algebraically closed field, and $X$ is a general set of $s$ points in $\PP^N_\kk$. Let $I$ be the defining ideal of $X$. If $m\in \NN$ and $w$ is the least integer such that $$(s-1){m+N-1 \choose N} < {N+w \choose N }$$ then $\reg(I^{(m)}) \leqslant m+w$.
\end{lemma}

\section{Demailly's conjecture large general points}\label{section: Demailly's conjecture for general points} 
In this section, we improve the known result of Demailly's conjecture for general points. It has been proved in \cite[Theorem 2.9]{bghndemailly} that for a fixed $m$, if $I$ is the ideal defining at least $(2m+2)^N$ general points in $\PP^N$ then 
$$ \ahat(I)\geqslant \dfrac{\alpha(I^{(m)})+N-1}{m+N-1}.$$ 
We start with the following binomial inequality inspired by \cite[Theorem 1]{chang2020demailly}.
\begin{lemma}\label{lemma: combinatorial inequality for general}
	Let  $k\geqslant 5$ and  $k^N\leqslant s < (k+1)^N$.	
	If $\displaystyle \sqrt[N]{s}=k+\epsilon$, where $0 \leqslant \epsilon <1$, and  $k \geqslant \dfrac{(1+\epsilon)2m}{N-1}+1+\epsilon$ then $${(k-1)(m+N-1)+N-1 \choose N}\geqslant s^N{m+N-1\choose N}.$$
\end{lemma}
\begin{proof}
	Now$${(k-1)(m+N-1)+N-1 \choose N} \geqslant s^N{m+N-1\choose N}$$ is equivalent to  
	$$ \left( (k-1)(m+N-1)+N-1 \right) \dots (k-1)(m+N-1)>(k+\epsilon)^N(m+N-1)\dots m.$$
	It is enough to show that, for $i=0, \dots, \left\lfloor \frac{N-1}{2} \right\rfloor$, \begin{equation}\label{eqn: first equation} \left((k-1)(m+N-1) +N-1-i \right)\left((k-1)(m+N-1)+i \right) \geqslant  (k+\epsilon)^2(m+N-1-i)(m+i)  .\end{equation}	
	Now LHS (\ref{eqn: first equation})-RHS(\ref{eqn: first equation})= $Am^2+Bm+C_i$, 
	where
	\begin{align*}
		A&=-2k\epsilon-\epsilon^2-2k+1=-(1+\epsilon)(2k-1+\epsilon)<0,\\
		B&= k^2N-2k\epsilon N-\epsilon^2N-k^2+2k\epsilon+\epsilon^2-3kN+3k+N-1=(N-1)(k^2-k(2\epsilon+3)+(1-\epsilon^2) )>0, \\
		C_i&=((k+\epsilon)^2-1)i^2-(N-1)((k+\epsilon)^2 -1)i+(N-1)^2(k-1)k.
	\end{align*}
	Consider $C_i$ as a quadratic in $i$. Since $\text{coeff}(i^2)>0$, the parabola is concave up. Its minimum is obtained at $i=-\dfrac{\text{coeff(i)}}{2\text{coeff}(i^2)}=\dfrac{N-1}{2}$ and the value is $C=\dfrac{(N-1)^2}{4}[(2k-1)^2-(k+\epsilon)^2]>0$.  \\
	To prove equation \ref{eqn: first equation}, it is enough to show that \begin{equation}\label{eqn: quadtatic Am^2+Bm+C}
		Am^2+Bm+C\geqslant 0, \text{ under the given condition}.
	\end{equation}
	
	Since $A<0$ and $C>0$, the parabola $Am^2+Bm+C$ is concave down. Again $B^2-4AC=(N-1)^2(k^2+k\epsilon)^2 \geqslant 0$. Thus, each coordinate of the vertex is positive. Thus, it is enough to show that the parabola $Am^2+Bm+C$ is positive when $0\leqslant m \leqslant \text{largest root of the parabola}$. 
	The roots of the parabola are $$\dfrac{-B\pm \sqrt{B^2-4AC}}{2A}=\dfrac{-\left((N-1)(k^2-k(2\epsilon+3)+(1-\epsilon^2)\right)\pm \sqrt{(N-1)^2(k^2+k\epsilon)^2 }}{-(1+\epsilon)(2k-1+\epsilon)}. $$ 
	So the biggest root is $\dfrac{(N-1)(k-1-\epsilon)}{1+\epsilon}>0.$

	The smallest root is $-\dfrac{(N-1)(3k-1+\epsilon)}{(2k-1+ \epsilon)}<0$. 
	Thus $Am^2+Bm+C \geqslant 0$ whenever $0 \leqslant m \leqslant \dfrac{(N-1)(k-1-\epsilon)}{1+\epsilon} $. 
	Hence  (\ref{eqn: first equation}) holds whenever $ m \geqslant  \dfrac{(N-1)(k-1-\epsilon)}{1+\epsilon} $, or, equivalently, $k \geqslant \dfrac{(1+\epsilon)(2m)}{N-1}+1+\epsilon$.  
\end{proof}
\begin{lemma}\label{lemma: strong containment Deamailly}
	Let $I$ denote the ideal defining $s$ generic points in $\PP^N$. Let $k^N \leqslant s < (k+1)^N$ and $\mm$ denotes the maximal homogeneous ideal in $\kk[\PP^N]$. If $k \geqslant \dfrac{(1+\epsilon)2m}{N-1}+1+\epsilon,$ where $\sqrt[N]{s}=k+\epsilon$, and $0\leqslant \epsilon <1$, then we have 
	$$I^{(r(m+N-1)-N+1)} \subset \mm^{r(N-1)}(I^{(m)})^r, \text{ for } r \gg 0.  $$
\end{lemma}

\begin{proof} 
	The proof has the same structure as \cite[Lemma 2.8]{bghndemailly}, so we briefly outline it. 
	We need to prove that, for $r\gg0$, $$ \alpha(I^{(r(m+N-1)-N+1)}) \geqslant r\reg(I^{(m)})+r(N-1). $$
	Now from \cite[Theorem 2.4]{TrungValla}  and Lemma \ref{lemma: combinatorial inequality for general} we have 
	\begin{align*}
		\reg(I^{(m)})<(k-1)(m+N-1)-1+m\\
		\reg(I^{(m)}) <k(m+N-1) -(N-1)-1\\
		\reg(I^{(m)})+N-1 <k(m+N-1)-\dfrac{k(N-1)}{r}, \text{ for } r\gg0\\
		r(\reg(I^{(m)})+N-1)<k(r(m+N-1) -N+1) \leqslant\alpha(I^{(r(m+N-1)-N+1)}). 
	\end{align*} 
\end{proof}
Finally, we get Demailly's conjecture as follows. 
\begin{theorem}\label{theorm: Demailly general}
	Let $\kk$ be an algebraically closed field and $N\geqslant 3$. Fix an integer $m\geqslant 1$. If $I$ denotes the ideal defining a set of at least $\left( \dfrac{(1+\epsilon)2m}{N-1}+1+\epsilon \right)^N $ general points in $\PP^N_\kk$, where $0\leqslant \epsilon <1$, then 
	
	$$\ahat(I)\geqslant\dfrac{\alpha(I^{(m)})+N-1}{m+N-1}.$$
\end{theorem}
\begin{proof}
	The proof follows from Lemma \ref{lemma: strong containment Deamailly} and it is similar as \cite[Theorem 2.9]{bghndemailly}.  
\end{proof}
\begin{remark}\label{remark: Demailly general what has been improved}
	Theorem \ref{theorm: Demailly general} improves the restriction on the number of points. For example, from \cite{bghndemailly}  Demailly's conjecture for $m=4$ is known for ideals defining at least $8^N$ many general points in $\PP^N, N\geqslant 9$. Theorem \ref{theorm: Demailly general} improves that for at least $ 4^N, N\geqslant 9$. In general, Theorem \ref{theorm: Demailly general} tells that for any $m$, ideal $I$ defining  at least $4^N$ general points in $\PP^N$, where $N\geqslant 2m+1$,  satisfies $\ahat(I)\geqslant \dfrac{\alpha(I^{(m)}+N-1)}{m+N-1}$.  Previously, it was known to be true for at least $(2m)^N $ general points \cite[Theorem 2.9]{bghndemailly}. 
\end{remark}

\section{Demailly's conjecture for $m=3$}\label{section: compairing waldschmidt constant with the third symbolic powers}  
In this section, we prove Demailly's conjecture for $m=3$, i.e., concerning the initial degree of the third symbolic power of the ideal defining points. We recall the statement of the conjecture. 

\vspace{0.5em}

\noindent\textbf{Conjecture \ref{conjecture: Demailly3}.}
If $I$ is the defining ideal of $s$ reduced points in $\PP^N$ then $$\ahat(I)\geqslant \dfrac{\alpha(I^{(3)})+N-1}{N+2}. $$

We start with the following corollary, a direct application of Theorem \ref{theorm: Demailly general}. 
\begin{corollary}\label{corollary: Demailly m=3 from demailly general}
	Let $I$ denote the ideal defining $s$ general points in $\PP^N$. If $s$ is any of the following: 
	\begin{enumerate}[label=(\roman*)]
		\item $s\geqslant 6^4$, where $N=4$
		\item $s \geqslant 5^N$, where $N=5$, or $6$
		\item $s\geqslant 4^N$, where $N \geqslant 7$.
	\end{enumerate}
	then the ideal $I$ satisfies Conjecture \ref{conjecture: Demailly3} $$\ahat(I)\geqslant \dfrac{\alpha(I^{(3)})+N-1}{N+2}.$$
\end{corollary}
\begin{proof}
	From Theorem \ref{theorm: Demailly general} we get for $m=3$, if 
	$s\geqslant \left( \dfrac{(1+\epsilon)6}{N-1}+1+\epsilon \right)^N$ then the ideal $I$ defining $s$ general points satisfies $$\ahat(I)\geqslant \dfrac{\alpha(I^{(3)})+N-1}{N+2}.$$ 
	Now, depending on the values of $N$, we get the lower bounds on $s$ as described in the hypothesis. 
\end{proof}
Now we investigate Conjecture \ref{conjecture: Demailly3} for very general points. We start with the following binomial inequalities. 
\begin{lemma} \label{lemma: combinatorial bound}
	We have the following inequalities.
	\begin{enumerate}
		\item $3^N  {N+2 \choose N}\leqslant {2N+6 \choose N}$ for $N\geqslant 8$.
		\item $2^N{N+2 \choose N} \leqslant {2N \choose N}$ for $N \geqslant 8$.
		\item $3^N  {N+2 \choose N}\leqslant {2N \choose N}$ for $N \geqslant 30$. 
	\end{enumerate}
\end{lemma}
\begin{proof}
	$(1)$ We prove by induction. For $N=8$, we can check that $3^8{10 \choose 8}< {22 \choose 8}$. Assume that the inequality holds for an arbitrary $N>8$ that is $3^N  {N+2 \choose N}\leqslant {2N+6 \choose N}$. Now 
	\begin{align*}
		3^{N+1}{N+3 \choose N+1}&=\dfrac{3(N+3)}{(N+1)} 3^N \cdot {N+2 \choose N}\\
		& \leqslant \dfrac{3(N+3)}{(N+1)} {2N+6 \choose N}; \text{ (by induction)}\\ 
		& =\dfrac{3(N+3)(N+1)}{(N+1)(2N+7)} {2N+6 \choose N}\\
		& =\dfrac{3(N+3)(N+7)}{(2N+7)(2N+8)}{2N+8 \choose N+1}\\
		& \leqslant {2N+8 \choose N+1}. 
	\end{align*}
	$(2)$ This follows the similar pattern as in $(1)$. \\
	$(3)$ Again one easy calculation shows that $3^{30}{32 \choose 20}< {60 \choose 30}$. 
	Assume that the inequality holds for an arbitrary $N >30$. The inductive steps are the same as in $(1)$, so we omit it.
\end{proof}

First, we prove Conjecture \ref{conjecture: Demailly3} for the following cases.

\begin{lemma} \label{thm: 2^N <= 3^N, 7 <= N<=29}
	Let $7 \leqslant N \leqslant 29$ and $I$ be the defining ideal of $s$ very general points in $\PP^N$. 
	If $2^N \leqslant s \leqslant 3^N$, then $$\ahat(I)\geqslant \dfrac{\alpha(I^{(3)})+N-1}{N+2}. $$
\end{lemma}
\begin{proof}
	Note that $$ \text{ if }{\ell+N\choose N}\leqslant  s{N+2\choose N} < {\ell+N+1\choose N} \text{ then } \alpha(I^{(3)})\leqslant \ell+1.  $$ 
	Thus it is enough to show that $$ \ahat({I}) \geqslant \dfrac{N+\ell}{N+2}, \text{ where } {\ell+N\choose N}\leqslant  s{N+2\choose N} < {\ell+N+1\choose N}. $$
	By Lemma \ref{lemma: combinatorial bound} part (1) and (2), we have $ N+1 \leqslant \ell \leqslant N+6$. For $\ell=N+1, \dots, N+4$, we have $\alpha(I)\geqslant 2 \geqslant \dfrac{N+\ell }{N+2}$. Now, with direct computation, that 
	$2^N (N+1) \leqslant \left \lceil  \dfrac {{2N+5 \choose N}}{{N+2 \choose N}}  \right\rceil $, when $7 \leqslant N \leqslant 29$. 
	Thus if $\ell=N+5$, or $N+6$, then using \cite[Theorem 3.2]{bisuinguyen2021chudnovsky} and \cite[Proposition B.1.1]{DUMNICKI2014471} we get
	$$\ahat(I)\geqslant \ahat \left( (N+1)2^N \right) \geqslant  2\cdot \dfrac{N+1}{N}>\dfrac{N+\ell}{N+2} .$$
\end{proof}
In the following theorem, we study a slightly weaker version of Conjecture \ref{conjecture: Demailly3}, and we prove that it is true for ideals defining every very general point. On the way, we prove Conjecture \ref{conjecture: Demailly3} for a certain number of points described in Remark \ref{remark: demailly for few cases}.  
\begin{theorem}\label{thm: demailly very general when m=3} 
	If $I$ is the defining ideal of $s$ very general points in $\PP^N$, then $$\ahat(I)\geqslant \dfrac{\alpha(I^{(3)})+N-2}{N+2}. $$
\end{theorem}

\begin{proof} 
	If $3 \leqslant N \leqslant 7$, then the result follows from Corollary \ref{corollary: demailly m=3 very general in lower PN}. 	
	Thus, we assume $N\geqslant 8$, and the results in Section \ref{section: computations in lower dimensions} are considered as base cases here. 
	It is enough to prove that $$\ahat(I)\geqslant \dfrac{N+\ell-1}{N+2}, \text{ where } {N+\ell \choose N} \leqslant s{N+2 \choose N} < {N+\ell+1 \choose N}.$$	
	Again it is already known when $s \leqslant N+2$ \cite[Corollary 5.7]{nagel2019interpolation} and $s\geqslant 3^N$ \cite[Theorem 4.8]{dumnicki2024waldschmidt}, so we prove for $N+3< s < 3^N$.  
	Note that $${N+4 \choose N} \leqslant (N+3){ N+2 \choose N } < {N+5 \choose N}, \text{ for } 5 \leqslant {N \leqslant 8}.$$ and
	$${N+3\choose N} \leqslant (N+3) {N+2 \choose N} < {N+4 \choose N}, \text{ for } N \geqslant 9.$$
	If  $\ell=3$ or $4$ and  ${N+\ell \choose N}  \leqslant s{N+2 \choose N} < {N+\ell \choose N}$ then $ \ahat(I) \geqslant \ahat(N+3) >\dfrac{N+2}{N}\geqslant \dfrac{N+\ell}{N+2}$. We can assume that $\ell \geqslant 5$. Again, by Lemma \ref{lemma: combinatorial bound}, we have $\ell \leqslant N-1$, for $N \geqslant 30$. Following, Lemma \ref{thm: 2^N <= 3^N, 7 <= N<=29}, only remaining cases are $s\leqslant 2^N$, when $ 7 \leqslant N \leqslant 29$. For those cases, from Lemma \ref{lemma: combinatorial bound} we have $\ell \leqslant N$.  
	Hence, we only need to prove for $5 \leqslant \ell \leqslant N-1$, and $N\geqslant 8$.   
	
	\noindent\textbf{Case 1:} Assume that $\ell \leqslant  \dfrac{5+\sqrt{8N+17}}{2} $.  
	First we show that $(N+1){N+\ell\choose N}\geqslant{N+2\choose N}{N+\ell-2\choose N}$. Now the above inequality is equivalent to $$2(N+\ell)(N+\ell-1) \geqslant(N+2)\ell(\ell-1), \text{ i.e., } \ell^2-5\ell-2(N-1)\leqslant 0.$$
	This is a parabola, and it is concave. Its zeros are $\dfrac{5\pm\sqrt{8N+17}}{2}$. Again, we have $5\leqslant \ell$. Thus $\ell^2-5\ell-2(N-1)\leqslant 0$ when $ \ell\leqslant \dfrac{5+\sqrt{8N+17}}{2}$.
	Now, in this case we have, if $s{N+2\choose N} \geqslant {N+\ell\choose N}$ then $(N+1){s}\geqslant {N+\ell-2\choose N}$. Thus from \cite[Theorem 3.2]{bn2024lower} we get $\ahat(I) \geqslant \dfrac{N+\ell-2}{N+1}>\dfrac{N+\ell-1}{N+2}$.   
	
	\vspace{1em}
	\noindent\textbf{Case 2:} Assume that $\ell \geqslant  \dfrac{5+\sqrt{8N+17}}{2}$ .\\ 
	Let $s$ satisfies ${N+\ell \choose N} \leqslant s{N+2 \choose N} < {N+\ell+1 \choose N}$. 
	We claim that $$\dfrac{{N+\ell \choose N}}{{N+2 \choose N}} \geqslant \dfrac{{N-1+\ell \choose N-1}}{{N+1 \choose N-1}}+\dfrac{{N-1+\ell-1 \choose N-1}}{{N+1 \choose N-1}}+1.$$ 
	The above inequality is equivalent to $$\dfrac{(N+\ell-2)!}{\ell!(N+2)!}\left(\ell^2-5\ell -(2N-2) \right) \geqslant \dfrac{1}{2}.$$
	Let $\ell_0=\left\lceil \dfrac{5+\sqrt{8N+17}}{2} \right\rceil+1$. 
	Now 
	$$\resizebox{1.05\hsize}{!}{$\dfrac{(N+\ell_0-2)!}{\ell_0!(N+2)!}\left(\ell_0^2-(2N+1)\ell_0 -(2N-2) \right)=\dfrac{(N+\ell_0-2)(N+\ell_0-3)\dots (N+3)(\ell_0^2-(2N+1)\ell_0 -(2N-2))}{\ell_0(\ell_0-1)\dots 3\cdot 2}$} $$
	$$=\left(\dfrac{(N+\ell_0-2)  }{\ell_0}\right)\cdot\left(\dfrac{(N+\ell_0-3) }{(\ell_0-1)}\right)\dots \left(\dfrac{(N+4)}{4\cdot 5}\right)\cdot \left(\dfrac{(N+3)}{3\cdot 6 }\right)\cdot\left(\ell_0^2-5\ell_0 -(2N-2) \right) \cdot \dfrac{1}{2}.$$ 
	Except for the last term, each term in the product is at least $1$, for $N\geqslant 16$, and $\ell_0^2-5\ell_0 -(2N-2)>0$ thus $\dfrac{(N+\ell_0-2)!}{\ell_0!(N+2)!}\left(\ell_0^2-5\ell_0 -(2N-2) \right)\geqslant \dfrac{1}{2}$. When $8 \leqslant N \leqslant 15$, then simple calculations will give that $ \dfrac{(N+\ell_0-2)!}{\ell_0!(N+2)!}\left(\ell_0^2-5\ell_0 -(2N-2) \right)>\dfrac{1}{2}$.
	
	Whenever, $\ell> \dfrac{5+\sqrt{8N+17}}{2}$ and ${N+\ell \choose N}\leqslant s{N+2 \choose N}$ then 
	$$ \dfrac{{N+\ell \choose N}}{{N+2 \choose N}} \geqslant \dfrac{{N-1+\ell \choose N-1}}{{N+1 \choose N-1}}+\dfrac{{N-1+\ell-1 \choose N-1}}{{N+1 \choose N-1}}+1.$$ 
	This implies, 
	$$s \geqslant \left \lceil \dfrac{{N+\ell \choose N}}{{N+2 \choose N}} \right\rceil \geqslant \left\lceil \dfrac{{N-1+\ell \choose N-1}}{{N+1 \choose N-1}}+\dfrac{{N-1+\ell-1 \choose N-1}}{{N+1 \choose N-1}} \right\rceil +1 \geqslant \left\lceil \dfrac{{N-1+\ell \choose N-1}}{{N+1 \choose N-1}} \right\rceil + \left\lceil \dfrac{{N-1+\ell-1 \choose N-1}}{{N+1 \choose N-1}} \right\rceil$$ 
	Let $r_1 \geqslant \left\lceil \dfrac{{N-1+\ell \choose N-1}}{{N+1 \choose N-1}} \right\rceil$, and $r_2 \geqslant \left\lceil \dfrac{{N-1+\ell-1 \choose N-1}}{{N+1 \choose N-1}} \right\rceil$ be very general points in $\PP^{N-1}$. 
	Then by induction hypothesis: $$\ahat(\PP^{N-1}, r_1)\geqslant \dfrac{N-1+\ell }{N+1}=a_1, \text{ and } \ahat(\PP^{N-1}, r_2)\geqslant \dfrac{N-2+\ell }{N+1}=a_2.$$
	Here $a_1>1$ and $a_2 \leqslant 2$, since $5 \leqslant \ell \leqslant N-1$. Therefore using Lemma \ref{lemma: Waldschmidtdecomp} we get:  
	$$\ahat(\PP^N,s) \geqslant \left(1-\dfrac{N+1}{N-1+\ell}\right) \dfrac{N-2+\ell}{N+1}+1 \geqslant \dfrac{N+\ell}{N+2}.$$
	The last inequality equals $\ell^2-5\ell+6 \geqslant 0$, which holds as $\ell\geqslant 5$.
\end{proof}
\begin{remark}\label{remark: demailly for few cases}
	Notice when $\left\lceil  \dfrac{{N+\ell\choose N} }{{N+2\choose N} }\right\rceil \leqslant s \leqslant 3^N $, where $\ell\geqslant \dfrac{5+\sqrt{8N+17}}{2},$ then Theorem \ref{thm: demailly very general when m=3} shows that the defining ideal $I$ of $s$ very general points in $\PP^N$ satisfies Demailly's conjecture for $m=3$ that is, $\ahat(I)\geqslant \dfrac{\alpha(I^{(3)})+N-1}{N+2}$. 
	
\end{remark}

\section{Demailly's conjecture for $m=3$ in lower dimensional projective spaces }\label{section: computations in lower dimensions} 
Now we prove Conjecture \ref{conjecture: Demailly3} for very general points $\PP^N$, where $3 \leqslant N\leqslant 7$.  First, we prove that Conjecture \ref{conjecture: Demailly3} holds for at least $6$ and, at most $3^3 $ very general points in $\PP^3$.
\begin{theorem} \label{Computations in P^3}
	If $I$ denotes the ideal defining $s$ very general points in $\PP^3$ and $6 \leqslant s \leqslant 3^3$ then 
	$$\ahat(I)\geqslant \dfrac{\alpha(I^{(3)}) + 2}{5}.$$
\end{theorem}
\begin{proof}
	It is enough to prove that $$\ahat(I)\geqslant \dfrac{3+\ell}{3+2}=\dfrac{3+\ell}{5}, \text{ where } {3+\ell \choose 3} < s{3+2 \choose 3} \leqslant {3+\ell+1\choose 3}.$$
	Since ${3+5 \choose 3} < 10\cdot 6 < {3+6 \choose 3}$ and ${3+9 \choose 3} < 10\cdot 3^3 < {3+9+1\choose 3}$, then $5 \leqslant \ell \leqslant 9$.
	\begin{enumerate}[label=Case (\roman*).] 
		\item When $\ell=5,s=6,7,8$. By \cite[Proposition 11]{dumnicki2015containments}, we have $\ahat(I)\geqslant\ahat(6)\geqslant \dfrac{12}{7} > \dfrac{3+5}{5}.$ 
		\item When $\ell=6,9 \leqslant s \leqslant 12$. By \cite[Proposition 11]{dumnicki2015containments}, we have $\ahat(I)\geqslant\ahat(2^3)\geqslant 2 \geqslant \dfrac{3+6}{5}$.
		\item When $\ell=7,13 \leqslant s \leqslant 16$. By \cite[Proposition 11]{dumnicki2015containments}, we have $$\ahat(I)\geqslant\ahat(12)\geqslant \dfrac{126}{57} \cong 2.21 > \dfrac{3+7}{5}.$$
		\item When $\ell=8,17 \leqslant s \leqslant 22$. By \cite[Proposition 11]{dumnicki2015containments}, we have $\ahat(I)\geqslant\ahat(17)\geqslant \dfrac{5}{2}> \dfrac{3+8}{5}$.
		\item When $\ell=9,23 \leqslant s \leqslant 3^3$. By \cite[Proposition 11]{dumnicki2015containments}, we have $\ahat(I)\geqslant\ahat(21)\geqslant \dfrac{8}{3} > \dfrac{3+9}{5}$.
		\end{enumerate}
\end{proof}
Now we  prove that, Conjecture \ref{conjecture: Demailly3} is satisfied in  $\PP^N$, where $4\leqslant N \leqslant 7$.
\begin{theorem} \label{Computations in P^4}
	If $I$ denotes the ideal defining $s$ very general points in $\PP^4$ and $7 \leqslant s \leqslant 3^4$ then 
	$$\ahat(I)\geqslant \dfrac{\alpha(I^{(3)}) + 3}{6}.$$
\end{theorem}
\begin{proof}
	It is enough to prove that $$\ahat(I)\geqslant \dfrac{4+\ell}{4+2}=\dfrac{4+\ell}{6}, \text{ where } {4+\ell \choose 4} < s{4+2 \choose 4} \leqslant {4+\ell+1\choose 4}.$$
	Since ${4+4 \choose 4} < 15\cdot 4 < {4+5 \choose 4}$ and ${4+10 \choose 4} < 15\cdot 3^4 < {4+11 \choose 4}$, then $4 \leqslant \ell \leqslant 10$.
	\begin{enumerate}[label=Case (\roman*).] 
		\item When $\ell=4,5 \text{ then }7\leqslant s \leqslant 14$. By \cite[Proposition B.1.1]{DUMNICKI2014471}, we have $$\ahat(I)\geqslant\ahat(\PP^4,7)=\ahat(\PP^4,6)\geqslant \dfrac{6}{4} \geqslant \dfrac{4+5}{6}.$$
		\item When $\ell=6,7 \text{ then }15 \leqslant s \leqslant 33$. Now by Theorem \ref{Computations in P^3} \[r_1= 9, a_1= 2; r_2=6, a_2 = \dfrac{12}{7}.\]
		Then $ s\geqslant r_1+r_2=15$. Now $1 < a_1 \leqslant 2$ and $a_2 \leqslant 2$ thus by Lemma \ref{lemma: Waldschmidtdecomp} we get 
		$$\ahat(I)\geqslant\ahat(\PP^4, 15)\geqslant (1-\dfrac{1}{a_1})a_2+1=(1-\dfrac{1}{2})\cdot \dfrac{12}{7}+1= 1.85>\dfrac{4+7}{6}\cong 1.83.$$
		\item When $\ell=8,34 \leqslant s \leqslant 47$.  We have by \cite[Lemma 2.5]{bn2024lower} $$\ahat(I)\geqslant\ahat(\PP^4,34)\geqslant\ahat(\PP^4,2^4)= 2 \geqslant \dfrac{4+8}{6}.$$
		\item When $\ell=9,10, \text{ then }48 \leqslant s \leqslant 3^4$. Now by Theorem \ref{Computations in P^3} $$r_1= 23, a_1 = \dfrac{8}{3}; r_2= 17, a_2 = \dfrac{5}{2} \text{ and } r_3= 6, a_3 = \dfrac{12}{7}. $$
		Now $s\geqslant r_1+r_2+r_3=46$. Now $2 < a_1,a_2 \leqslant 3$ and $a_3 \leqslant 3$ thus by Lemma \ref{lemma: Waldschmidtdecomp} and \cite[Proposition B.1.1]{DUMNICKI2014471} we get 
		$$\ahat(I)\geqslant\ahat(\PP^4, 46)\geqslant (1-\dfrac{1}{a_1}-\dfrac{1}{a_2})a_3+2=(1-\dfrac{3}{8}-\dfrac{2}{5})\dfrac{12}{7}+2\cong 2.38>\dfrac{4+10}{6}\cong 2.33.$$

	\end{enumerate}
\end{proof}

\begin{theorem} \label{Computations in P^5}
	If $I$ denotes the ideal defining $s$ very general points in $\PP^5$ and $8 \leqslant s \leqslant 3^5$ then  
	$$\ahat(I)\geqslant \dfrac{\alpha(I^{(3)})+4}{7}.$$
	
\end{theorem}
\begin{proof} 
	It is enough to show that 	$$\ahat(I)\geqslant \dfrac{5+\ell}{5+2}=\dfrac{5+\ell}{7}, \text{ where } {5+\ell \choose 5} < s{5+2 \choose 5} \leqslant {5+\ell+1\choose 5}.$$ 
	Since ${5+4 \choose 5} < 21\cdot 8 < {5+5 \choose 5}$ and ${5+11 \choose 5}< 21\cdot 3^5 < {5+12 \choose 5}$, then $4 \leqslant \ell \leqslant 11$.
	\begin{enumerate}[label=Case (\roman*).] 
		\item When $\ell=4,8\leqslant s\leqslant 12$. By \cite[Proposition B.1.1]{DUMNICKI2014471}, we have $$\ahat(I)\geqslant\ahat(\PP^5,8)\geqslant \dfrac{24}{17} > \dfrac{5+4}{7}.$$
		\item When $\ell=5, 13\leqslant s \leqslant 22$.  Now by Theorem \ref{Computations in P^4} and \cite[Proposition B.1.1]{DUMNICKI2014471} $$r_1= r_2=6, a_1=a_2 = \dfrac{3}{2}.$$
		Then $ s\geqslant r_1+r_2=12$. Now $1 < a_1 \leqslant 2$ and $a_2 \leqslant 2$ thus by Lemma \ref{lemma: Waldschmidtdecomp} we get 
		$$\ahat(I)\geqslant\ahat(\PP^5, 13)\geqslant (1-\dfrac{1}{a_1})a_2+1=(1-\dfrac{2}{3})\cdot \dfrac{3}{2}+1= 1.5>\dfrac{5+5}{7}\cong 1.43.$$
		\item When $\ell=6, \text{ then } 23\leqslant s \leqslant37$.  Now by Theorem \ref{Computations in P^4} 
		$$r_1= 15,a_1 = 1.85 \text{ and }r_2= 7,a_2 = 1.5.$$
		Now $s\geqslant r_1+r_2=22$. Now $1 < a_1 \leqslant 2$ and $a_2 \leqslant 2$ thus by Lemma \ref{lemma: Waldschmidtdecomp} we get 
		$$\ahat(I)\geqslant\ahat(\PP^5, 23)\geqslant (1-\dfrac{1}{a_1})a_2+1=(1-\dfrac{1}{1.85})\cdot 1.5+1 \cong 1.689>\dfrac{5+6}{7}\cong 1.57.$$
		\item When $\ell=7,8,9$, then $38\leqslant s \leqslant 143$ and thus $\ahat(I)\geqslant \ahat(\PP^5,2^5)= 2 \geqslant \dfrac{5+9}{7}$.
		\item When $\ell=10, 11$, then  $144 \leqslant s \leqslant 3^5$. 
		Now by Theorem \ref{Computations in P^4}  $$ r_1= 48, a_1 = 2.38; r_2= 48, a_2 =2.38 \text{ and } r_3= 34, a_3 = 2.$$ 
		Now $s> r_1+r_2+r_3=130$. Now $1 < a_1, a_2 <3$ and $a_3<3$ thus by Lemma \ref{lemma: Waldschmidtdecomp} we get 
		$$\ahat(I)\geqslant\ahat(\PP^5, 130)\geqslant (1-\dfrac{1}{a_1}-\dfrac{1}{a_2})a_3+2=(1-\dfrac{1}{2.38}-\dfrac{1}{2.38})\cdot2+2\cong 2.31>\dfrac{5+11}{7}\cong 2.29.$$
	
	\end{enumerate}
\end{proof}

\begin{theorem} \label{Computations in P^6}
	If $I$ denotes the ideal defining $s$ very general points in $\PP^6$ and $9 \leqslant s \leqslant 3^6$ then 
	$$\ahat(I)\geqslant \dfrac{\alpha(I^{(3)})+5}{8}.$$
	
\end{theorem}
\begin{proof} 
	It is enough to prove that 	$$\ahat(I)\geqslant \dfrac{6+\ell}{6+2}=\dfrac{6+\ell}{8}, \text{ where } {6+\ell \choose 6} < s{6+2 \choose 6} \leqslant {6+\ell+1\choose 6}.$$
	Since ${6+4 \choose 6} < 28\cdot 9 < {6+5 \choose 6}$ and ${6+12 \choose 6}< 28\cdot 3^6 < {6+13 \choose 6}$, then $4 \leqslant \ell \leqslant 12$.
	\begin{enumerate}[label=Case (\roman*).] 
		\item When $\ell=4,9\leqslant s\leqslant 16$. By \cite[Proposition B.1.1]{DUMNICKI2014471}, we have $$\ahat(I)\geqslant\ahat(\PP^6,9)\geqslant \dfrac{63}{47} > \dfrac{6+4}{8}.$$
		\item When $\ell=5, 17\leqslant s \leqslant 33$.  Now by Theorem \ref{Computations in P^5} we have $$r_1= 8,r_2=8, a_1=a_2 = 1.4117. $$
		Now $s\geqslant r_1+r_2=16$. Now $1 < a_1 \leqslant 2$ and $a_2 \leqslant 2$ thus by Lemma \ref{lemma: Waldschmidtdecomp} we get 
		$$\ahat(I)\geqslant \ahat(\PP^6, 16) \geqslant (1-\dfrac{1}{a_1})a_2+1=(1-\dfrac{1}{1.4117})\cdot 1.4117+1= 1.4117>\dfrac{6+5}{8}.$$
		\item When $\ell=6, 34\leqslant s \leqslant 61$.  Now by Theorem \ref{Computations in P^5} we get
		$$ r_1= 23,a_1 = 1.689 \text{ and } r_2=8, a_2 = 1.4117. $$
		Now $s\geqslant r_1+r_2=31$, $1 < a_1 \leqslant 2$ and $a_2 \leqslant 2$ thus by Lemma \ref{lemma: Waldschmidtdecomp} we get 
		$$\ahat(I)\geqslant\ahat(\PP^6, 31)\geqslant (1-\dfrac{1}{a_1})a_2+1=(1-\dfrac{1}{1.689})\cdot 1.4117+1 \cong 1.5758>\dfrac{6+6}{8}.$$
		\item When $\ell=7, 62\leqslant s \leqslant107$.  Now by Theorem \ref{Computations in P^5} we get
		$$ r_1= 32,a_1 = 2 \text{ and } r_2=23, a_2 = 1.689. $$
		Now $s\geqslant r_1+r_2=55$, $1 < a_1 \leqslant 2$ and $a_2 \leqslant 2$ thus by Lemma \ref{lemma: Waldschmidtdecomp} we get 
		$$\ahat(I)\geqslant \ahat(\PP^6, 55)\geqslant (1-\dfrac{1}{a_1})a_2+1=(1-\dfrac{1}{2})\cdot 1.689+1 \cong 1.8445>\dfrac{6+7}{8}.$$
		\item When $\ell=8,\ldots,10$ we have $108\leqslant s \leqslant442$.  Now by \cite[Lemma 2.5]{bn2024lower} we get
		$$\ahat(\PP^6,108)\geqslant \ahat(\PP^6,2^6)=2\geqslant \dfrac{6+10}{8}.$$
		\item When $\ell=11,12$, then  $443 \leqslant s \leqslant 3^6$. 
		$$ r_1= 144,a_1 = 2.31; r_2=144, a_2= 2.31\text{ and } r_3=32, a_3 = 2 .$$
		Now $s\geqslant r_1+r_2+r_3=320$, $2 < a_1,a_2 \leqslant 3$ and $a_3 \leqslant 3$ thus by Lemma \ref{lemma: Waldschmidtdecomp} we get 
		$$\ahat(I)\geqslant\ahat(\PP^6,320)\geqslant (1-\dfrac{1}{a_1}-\dfrac{1}{a_2})a_3+2=(1-\dfrac{1}{2.31}-\dfrac{1}{2.31})\cdot 2+2 \cong 2.268>\dfrac{6+12}{8}.$$ 
	\end{enumerate}
\end{proof}

\begin{theorem} \label{Computations in P^7}
	If $I$ denotes the ideal defining $s$ very general points in $\PP^7$ and $10 \leqslant s \leqslant 3^7$ then 
	$$\ahat(I) \geqslant \dfrac{\alpha(I^{(3)})+6}{9}.$$
\end{theorem}
\begin{proof} 
	It is enough to show that $$\ahat(I)\geqslant \dfrac{7+\ell}{7+2}=\dfrac{7+\ell}{9}, \text{ where } {7+\ell \choose 7} < s{7+2 \choose 7} \leqslant {7+\ell+1\choose 7}.$$ 
	Since ${7+4 \choose 7} < 36\cdot 10 < {7+5 \choose 7}$ and ${7+13 \choose 7}< 36\cdot 3^7 < {7+14 \choose 7}$, then $4 \leqslant \ell \leqslant 13$.
	\begin{enumerate}[label=Case (\roman*).] 
		\item When $\ell=4,10\leqslant s\leqslant 22$. By \cite[Proposition B.1.1]{DUMNICKI2014471}, we have $$\ahat(I)\geqslant\ahat(\PP^7, 10)\geqslant 1.2903 > \dfrac{7+4}{9}.$$
		\item When $\ell=5, 23\leqslant s \leqslant47$.  Now by Theorem \ref{Computations in P^6} we have $$r_1=r_2=9, a_1=a_2 = \dfrac{63}{47}. $$
		Now $s\geqslant r_1+r_2=18$. Now $1 < a_1 \leqslant 2$ and $a_2 \leqslant 2$ thus by Lemma \ref{lemma: Waldschmidtdecomp} we get 
		$$\ahat(I)\geqslant\ahat(\PP^7, 18)\geqslant  (1-\dfrac{1}{a_1})a_2+1=(1-\dfrac{47}{63})\cdot \dfrac{63}{47}+1= \dfrac{63}{47}>\dfrac{7+5}{9}.$$
		\item When $\ell=6, 48\leqslant s \leqslant95$.  Now by Theorem \ref{Computations in P^6} we get
		$$ r_1= 34,a_1 = 1.5758 \text{ and } r_2=9, a_2 = \dfrac{63}{47}. $$
		Now $s\geqslant r_1+r_2=43$, $1 < a_1 \leqslant 2$ and $a_2 \leqslant 2$ thus by Lemma \ref{lemma: Waldschmidtdecomp} we get 
		$$\ahat(I)\geqslant\ahat(\PP^7, 43)\geqslant (1-\dfrac{1}{a_1})a_2+1=(1-\dfrac{1}{1.5758})\cdot \dfrac{63}{47}+1 \cong 1.4897>\dfrac{7+6}{9}.$$
		\item When $\ell=7,8, $\text{ then }$ 96\leqslant s \leqslant317$.  Now by Theorem \ref{Computations in P^6} 
		$$ r_1= 62,a_1 = 1.8445 \text{ and } r_2=34, a_2 = 1.5758. $$
		Now $s\geqslant r_1+r_2=96$, $1 < a_1 \leqslant 2$ and $a_2 \leqslant 2$ thus by Lemma \ref{lemma: Waldschmidtdecomp} we get 
		$$\ahat(I)\geqslant\ahat(\PP^7,96)\geqslant  (1-\dfrac{1}{a_1})a_2+1=(1-\dfrac{1}{1.8445})\cdot 1.5758+1 \cong 1.7214>\dfrac{7+8}{9}.$$
		\item When $\ell=9,\ldots,11$ we have $318\leqslant s \leqslant1399$.  Now  by \cite[Lemma 2.5]{bn2024lower}
		$$\ahat(\PP^7,318)\geqslant \ahat(\PP^7,2^7)=2\geqslant \dfrac{7+11}{9}.$$
		\item When $\ell=12,13$, then  $1400 \leqslant s \leqslant 3^7$. 
		Now by \cite[Corollary 3.5]{bisuinguyen2021chudnovsky} and \cite[Proposition B.1.1]{DUMNICKI2014471}, \vspace{-.4em} $$\ahat(\PP^7,1400)\geqslant \ahat(\PP^7, 10\cdot(2^7)^1)\geqslant 2\cdot\ahat(\PP^7,10)\geqslant2\cdot1.2903 \cong 2.5806>\dfrac{7+13}{9}.$$ \vspace{-2em}
	
	\end{enumerate}
\end{proof}
\vspace{0em}
Demailly's conjecture is known for at most $N+2$ points \cite[Corollary 5.7]{nagel2019interpolation} and at least $3^N $  \cite[Theorem 4.8]{dumnicki2024waldschmidt} very general  points in $\PP^N$.  Thus from Theorem \ref{Computations in P^3}, \ref{Computations in P^4}, \ref{Computations in P^5}, \ref{Computations in P^6} and \ref{Computations in P^7} we get the following corollary. 
\begin{corollary}\label{corollary: demailly m=3 very general in lower PN}
	If $I$ is the defining ideal of $s$ very general points in $\PP^N$, where $3 \leqslant N \leqslant 7$, then $I$ satisfies Demailly's conjecture for $m=3$, that is $$
	\ahat(I)\geqslant \dfrac{\alpha(I^{(3)})+N-1}{N+2}.$$
\end{corollary}

	\bibliographystyle{alpha}
	\bibliography{References}
	
\end{document}